\documentclass{amsart}

\newtheorem{theorem}{Theorem}[section]
\newtheorem{lemma}[theorem]{Lemma}

\theoremstyle{definition}

\theoremstyle{remark}

\numberwithin{equation}{section}

\newcommand{\bsx}{\boldsymbol{x}}
\newcommand{\NN}{\mathbb{N}}
\newcommand{\RR}{\mathbb{R}}
\newcommand{\Scal}{\mathcal{S}}
\newcommand{\Xcal}{\mathcal{X}}

\begin{document}

\title{The Sobol' sequence is not quasi-uniform in dimension 2}

\author{Takashi Goda}
\address{School of Engineering, The University of Tokyo, 7-3-1 Hongo, Bunkyu-ku, Tokyo 113-8656, Japan}
\email{goda@frcer.t.u-tokyo.ac.jp}
\thanks{The work of the author is supported by JSPS KAKENHI Grant Number 23K03210.}


\date{}

\dedicatory{}

\commby{}

\begin{abstract}
    Are common quasi-Monte Carlo sequences quasi-uniform? While this question remains widely open, in this short note, we prove that the two-dimensional Sobol' sequence is not quasi-uniform. This result partially answers an unsolved problem of Sobol' and Shukhman (2007) in a negative manner.
\end{abstract}

\maketitle

\section{Introduction}
For a dimension $d\geq 1$, let $\Xcal$ be a compact subset of $\RR^d$ with $\mathrm{vol}(\Xcal)>0.$ Given an infinite sequence of points $\Scal=(\bsx_i)_{i\geq 0}$ in $\Xcal$, we denote the first $n$ points of $\Scal$ by $P_n=(\bsx_i)_{0\leq i<n}$.
The fill distance and the separation radius of $P_n$ are defined, respectively, as follows:
\[ h(P_n):= \sup_{\bsx\in \Xcal}\min_{0\leq i< n}\|\bsx-\bsx_i\|, \]
and
\[ q(P_n):=\frac{1}{2}\min_{0\leq i<j< n}\|\bsx_i-\bsx_j\|,\]
where $\|\cdot\|$ denotes the Euclidean distance on $\RR^d$. The fill distance $h$ is sometimes referred to as the covering radius, dispersion, or minimax distance, while the separation radius $q$ is also referred to as the packing radius or maximin distance. A sequence $\Scal$ is said to be \emph{quasi-uniform} over $\Xcal$ if there exists a constant $C>1$, independent of $n$, such that
\[ 1\leq \frac{h(P_n)}{q(P_n)}\leq C\]
holds for all $n\geq 1$. Note that the ratio $h/q$ is termed the mesh ratio, providing a measure of how uniformly the points in $\Scal$ are distributed over $\Xcal$.

Obviously, the decay rates of $h$ and $q$ should coincide for $\Scal$ to be quasi-uniform. Here, both $h$ and $q$ are non-increasing functions with respect to $n$, and, for any sequence $\Scal$, $h$ is known to decay no faster than $O(n^{-1/d})$, see, for instance, \cite[Lemma~2.1]{PZ23}. On the other hand, $q$ can decay at an arbitrarily fast rate; it can even become $0$ for all $n\geq n_0$ if there exist two indices $i,j$ with $0\leq i<j<n_0$ such that $\bsx_i=\bsx_j$, which makes the mesh ratio unbounded anymore. Recently, it has been proven in \cite{PZ23} that the greedy packing algorithm can construct quasi-uniform sequences with $C\leq 2$ for general $\Xcal\subset \RR^d$, and that, for the case $\Xcal=[0,1]^d$, the resulting sequences are given explicitly in $d=1,2,4$ by placing the initial point $\bsx_0$ at the center of $[0,1]^d$.

Quasi-uniform sequences are useful in some applications such as kernel interpolation and Gaussian process regression, where the mesh ratio of sampling nodes plays a pivotal role in deriving upper bounds on the approximation error \cite{NWW06,SW06,T20,WBG21}. 

In this note, we focus on the $d$-dimensional unit cube $\Xcal=[0,1]^d$, where quasi-Monte Carlo (QMC) sequences, or low-discrepancy sequences, are a promising candidate for explicit construction of quasi-uniform sequences. In his seminal paper \cite{N88}, Niederreiter showed an upper bound on the fill distance $h$ of $O(n^{-1/d})$ for $\Scal$ being a $(t,d)$-sequence in base $b$, including Sobol', Faure and Niederetier sequences as special cases \cite{S67,F82,N88}. Consequently, in order to determine whether a $(t,d)$-sequence is quasi-uniform or not, it suffices to show whether a lower bound on the separation radius $q$ of order $\Omega(n^{-1/d})$ holds. In their study \cite{SS07}, Sobol' and Shukhman considered the Sobol' sequence, a widely-known $(t,d)$-sequence in base $2$, and proved a lower bound $q(P_n)\geq \sqrt{d}/(2n).$ Although this result is not sufficient to determine whether the Sobol' sequence is quasi-uniform for $d\geq 2$, they also conducted numerical experiments suggesting that $q\asymp n^{-1/d}$ holds, which would imply that the Sobol' sequence is quasi-uniform for all dimensions. If this is true, the Sobol' sequence would find a theoretical justification for its use as sampling nodes for kernel interpolation or Gaussian process regression, beyond its typical use in high-dimensional numerical integration. To the best of the author's knowledge, however, there has been no significant progress in this research direction since then. Indeed, it is stated in \cite{WBG21} that 
\begin{quote}
\textit{Since quasi-uniformity as defined above is not studied in QMC, it is unclear when common QMC point sets are quasi-uniform.}
\end{quote}

In this short note, with an aim to attract more attention to this problem, we provide a simple proof that, for the two-dimensional Sobol' sequence, there exist infinitely many $n$'s for which $q(P_n)= 1/(\sqrt{2}(n+1))$. This finding implies that the order of the lower bound proven by Sobol' and Shukhman cannot be improved for the case of $d=2$, concluding that the Sobol' sequence is not quasi-uniform in dimension 2. The question of whether higher-dimensional Sobol' sequences and other QMC sequences exhibit quasi-uniformity remains a topic for future research.

\section{Results}
For the definition of the Sobol' sequence in a general dimension $d\geq 1,$ we refer to \cite[Chapter~8]{DP10}. Here we only introduce its two-dimensional version. For a non-negative integer $n$, denote the binary expansion of $n$ by
\[ n=n_0+n_1 2+n_2 2^2+\cdots\quad \text{with $n_0,n_1,n_2,\ldots \in \{0,1\}$,}\]
which is indeed a finite expansion. Then the $n$-th point $\bsx_n=(x_{n,1},x_{n,2})\in [0,1]^2$ is given by
\[ x_{n,j}=\frac{x_{n,j,1}}{2}+\frac{x_{n,j,2}}{2^2}+\frac{x_{n,j,3}}{2^3}+\cdots,\]
with 
\[ \begin{pmatrix} x_{n,j,1} \\ x_{n,j,2} \\ x_{n,j,3} \\ \vdots \end{pmatrix} = C_j \begin{pmatrix} n_0 \\ n_1 \\ n_2 \\ \vdots \end{pmatrix} \bmod 2, \]
for $j=1,2$, where the modulo operation is applied componentwise to a vector and $C_j\in \{0,1\}^{\NN\times \NN}$ denotes the generating matrix for the $j$-th coordinate. Here, the generating matrices $C_1$ and $C_2$ are explicitly given by the identity matrix and the Pascal matrix over the two-element field, respectively, i.e., with $C_1=(c_{i,j}^{(1)})$ and $C_2=(c_{i,j}^{(2)})$, they are given by
\[ c_{i,j}^{(1)}=\begin{cases} 1 & \text{if $i=j$,} \\ 0 & \text{otherwise,}\end{cases} \quad \text{and}\quad c_{i,j}^{(2)}=\binom{j-1}{i-1}\bmod 2, \]
where we set $\binom{n}{m}=0$ if $n<m.$ It is obvious that $C_2$ is an upper triangular matrix.

To prove our main result, the following property of the Pascal matrix is crucial.
\begin{lemma}\label{lem:parity}
    Let $m=2^{w}-1$ for an integer $w\geq 1$. For any $1\leq i\leq m$, it holds that
    \[ \sum_{j=1}^{m}c_{i,j}^{(2)}\bmod 2 = 1. \]
\end{lemma}
\begin{proof}
    For any $m\geq 1$ and $1\leq i\leq m$, it follows from the hockey-stick identity that
    \[ \sum_{j=1}^{m}c_{i,j}^{(2)}\bmod 2 = \sum_{j=i}^{m}\binom{j-1}{i-1}\bmod 2= \binom{m}{i}\bmod 2.\]
    Denote the binary expansions of $m$ and $i$ by
    \[ m=m_0+m_12+\cdots+m_{w-1} 2^{w-1},\]
    and
    \[ i=i_0+i_12+\cdots+i_{w-1} 2^{w-1},\]
    respectively, for some $w\geq 1$ with $m_0,\ldots,m_{w-1},i_0,\ldots,i_{w-1}\in \{0,1\}$. Then, Lucas' theorem tells us 
    \[ \binom{m}{i}\bmod 2= \prod_{j=0}^{w-1}\binom{m_j}{i_j}\bmod 2. \]
    If $m$ is given by the form $2^{w}-1$, we have $m_0=m_1=\cdots=m_{w-1}=1$. Then, regarding each factor of the last product, we have $\binom{m_j}{i_j}=\binom{1}{i_j}=1$ regardless of whether $i_j=0$ or $i_j=1$. Therefore it is concluded that 
    \[ \sum_{j=1}^{m}c_{i,j}^{(2)}\bmod 2= \prod_{j=0}^{w-1}\binom{m_j}{i_j}\bmod 2= 1 \]
    holds for any $1\leq i\leq m$.
\end{proof}

This lemma says that, when $m=2^w-1$, every row in the upper-left $m\times m$ sub-matrix of $C_2$ contains an odd number of $1$'s in its elements. Using this parity property, we can prove our main result of this note. 

\begin{theorem}
    For the two-dimensional Sobol' sequence $\Scal$, we have
    \[ q(P_{n})= \frac{1}{\sqrt{2}(n+1)}, \]
    if $n=2^{m}-1$ where $m$ is given by the form $2^w-1$ with an integer $w>1$.
\end{theorem}

\begin{proof}
    Let us consider the last point $\bsx_{n-1}$ in $P_n$, which is the $(n-1)$-th point in the sequence $\Scal$ (given that the sequence starts with the ``zero-th'' point $\bsx_0$). Since we assume that $n$ is in the form $2^{m}-1$, the binary expansion of $n-1=2^{m}-2$ is given by $n-1=2+2^2+\cdots+2^{m-1}$, so that $(n-1)_0=0$, $(n-1)_1=\cdots=(n-1)_{m-1}=1$, and $(n-1)_{m}=(n-1)_{m+1}=\cdots=0$. By the definition of the two-dimensional Sobol' sequence, the first coordinate of $\bsx_{n-1}$ is given by
    \[ x_{n-1,1}=\frac{(n-1)_0}{2}+\frac{(n-1)_1}{2^2}+\cdots=\sum_{i=1}^{m-1}\frac{1}{2^{i+1}}=\frac{1}{2}-\frac{1}{2^m}.\]
    Regarding the second coordinate of $\bsx_{n-1}$, the observation that $C_2$ is an upper triangular matrix trivially leads to $x_{n-1,2,i}=0$ for all $i>m$. For $i=1$, as it holds from the definition of $C_2$ that $c_{1,j}^{(2)}=1$ for all $j\geq 1$ and $m=2^w-1$ is odd, we have
    \begin{align*}
    x_{n-1,2,1} & =c_{1,1}^{(2)}\cdot (n-1)_0+c_{1,2}^{(2)}\cdot (n-1)_1+\cdots+c_{1,m}^{(2)}\cdot (n-1)_{m-1} \bmod 2 \\
    & = m-1 \bmod 2 = 0.
    \end{align*}
    For any $1<i\leq m$, by using the observation that $C_2$ is an upper triangular matrix again, it follows from Lemma~\ref{lem:parity} that the number of $1$'s in $c_{i,i}^{(2)},c_{i,i+1}^{(2)},\ldots,c_{i,m}^{(2)}$ is odd, giving
    \begin{align*}
        x_{n-1,2,i} & = c_{i,i}^{(2)}\cdot (n-1)_{i-1}+\cdots+c_{i,m}^{(2)}\cdot (n-1)_{m-1}\bmod 2\\
        & = c_{i,i}^{(2)}+\cdots+c_{i,m}^{(2)}\bmod 2=1.
    \end{align*}
    Thus we have
    \[ x_{n-1,2}=\frac{x_{n-1,2,1}}{2}+\frac{x_{n-1,2,2}}{2^2}+\cdots = \sum_{i=2}^{m}\frac{1}{2^i}=\frac{1}{2}-\frac{1}{2^m}.\]
    Noticing that $\bsx_1=(1/2,1/2)$, the separation radius of $P_{n}$ is bounded above by
    \[ q(P_{n})=\frac{1}{2}\min_{0\leq i<j< n}\|\bsx_i-\bsx_j\|\leq \frac{1}{2}\|\bsx_1-\bsx_{n-1}\|=\frac{1}{2^{m+1/2}}=\frac{1}{\sqrt{2}(n+1)}. \]

    Let us consider the first $2^m=n+1$ points of the two-dimensional Sobol' sequence. Since each one-dimensional projection of these points is given by a permutation of $\{k/2^m\mid 0\leq k<2^m\}$, for any pair of indices $i,j$ with $0\leq i<j<2^m$, it holds that
    \[ \|\bsx_i-\bsx_j\| \geq \left(\left(\frac{1}{2^m}\right)^2+\left(\frac{1}{2^m}\right)^2\right)^{1/2}=\frac{1}{2^{m-1/2}}=\frac{\sqrt{2}}{n+1}.\]
    This leads to the matching lower bound on $q(P_{n})$ as
    \[ q(P_{n})\geq q(P_{2^m})=\frac{1}{2}\min_{0\leq i<j< 2^m}\|\bsx_i-\bsx_j\|\geq \frac{1}{\sqrt{2}(n+1)}. \]
    Thus the result follows.
\end{proof}

It is evident that our proof relies heavily on the specific nature of $C_2$ being the Pascal matrix. Extending our argument to determine whether the $d$-dimensional Sobol' sequence is quasi-uniform for $d\geq 3$ seems to be not so straightforward.

\bibliographystyle{amsplain}
\bibliography{ref.bib}

\providecommand{\bysame}{\leavevmode\hbox to3em{\hrulefill}\thinspace}
\providecommand{\MR}{\relax\ifhmode\unskip\space\fi MR }
\providecommand{\MRhref}[2]{%
  \href{http://www.ams.org/mathscinet-getitem?mr=#1}{#2}
}
\providecommand{\href}[2]{#2}
\begin{thebibliography}{10}

\bibitem{DP10}
J.~Dick and F.~Pillichshammer, \emph{{Digital nets and sequences: discrepancy theory and quasi-Monte Carlo integration}}, Cambridge University Press, 2010.

\bibitem{F82}
H.~Faure, \emph{Discrépance de suites associées à un système de numération (en dimension s)}, Acta Arith. \textbf{41} (1982), no.~4, 337--351 (in French).

\bibitem{NWW06}
F.~J. Narcowich, J.~D. Ward, and H.~Wendland, \emph{{Sobolev error estimates and a Bernstein inequality for scattered data interpolation via radial basis functions}}, Constr. Approx. \textbf{24} (2006), 175--186.

\bibitem{N88}
H.~Niederreiter, \emph{Low-discrepancy and low-dispersion sequences}, J. Number Theory \textbf{30} (1988), no.~1, 51--70.

\bibitem{PZ23}
L.~Pronzato and A.~Zhigljavsky, \emph{Quasi-uniform designs with optimal and near-optimal uniformity constant}, J. Approx. Theory \textbf{294} (2023), 105931.

\bibitem{SW06}
R.~Schaback and H.~Wendland, \emph{Kernel techniques: from machine learning to meshless methods}, Acta Numer. \textbf{15} (2006), 543--639.

\bibitem{S67}
I.~M. Sobol', \emph{On the distribution of points in a cube and the approximate evaluation of integrals}, Zh. Vychisl. Mat. Mat. Fiz. \textbf{7} (1967), no.~4, 784--802.

\bibitem{SS07}
I.~M. Sobol and B.~V. Shukhman, \emph{Quasi-random points keep their distance}, Math. Comput. Simulation \textbf{75} (2007), no.~3-4, 80--86.

\bibitem{T20}
A.~L. Teckentrup, \emph{{Convergence of Gaussian process regression with estimated hyper-parameters and applications in Bayesian inverse problems}}, SIAM/ASA J. Uncertain. Quantif. \textbf{8} (2020), no.~4, 1310--1337.

\bibitem{WBG21}
G.~Wynne, F.-X. Briol, and M.~Girolami, \emph{{Convergence guarantees for Gaussian process means with misspecified likelihoods and smoothness}}, J. Mach. Learn. Res. \textbf{22} (2021), no.~1, 5468--5507.

\end{thebibliography}

\end{document}